\documentclass[12pt,a4paper]{article}
\usepackage[latin2]{inputenc}
\usepackage{amsmath}
\usepackage{amsfonts}
\usepackage{amssymb}
\usepackage{graphicx}
\usepackage[left=2.00cm, right=2.00cm, top=2.50cm, bottom=2.50cm]{geometry}
\usepackage[T1]{fontenc}

\def\Q{{\mathbb Q}}
\def\Z{{\mathbb Z}}
\def\R{{\mathbb R}}

\newtheorem{lemma}{Lemma}
\newtheorem{theorem}[lemma]{Theorem}

\title{
Monogenity in totally real extensions \\of imaginary quadratic fields \\with
an application to simplest quartic fields
}
\author{
Istv\'{a}n Ga\'{a}l \\
{\small University of Debrecen, Mathematical Institute} \\
{\small H--4002 Debrecen Pf.400., Hungary,} \\
{\small e--mail: gaal.istvan@unideb.hu},
}

\begin{document}
\baselineskip=17pt

\maketitle
\thispagestyle{empty}

\renewcommand{\thefootnote}{\arabic{footnote}}
\setcounter{footnote}{0}

\vspace{0.5cm}

\begin{center}
{\it Dedicated to Professor Michael Pohst on his 75th birthday.}
\end{center}

\vspace{0.5cm}


\begin{abstract}
We describe an efficient algorithm to calculate generators of power integral bases
in composites of totally real fields with imaginary quadratic fields.
We show that the calculation can be reduced to solving index form equations 
in the original totally real fields.

We illustrate our method by investigating monogenity in 
the infinite parametric family of imaginary quadratic
extensions of the simplest quartic fields.
\end{abstract}

\newpage

\section{Introduction}


There is an extensive literature of monogenity of number fields and power integral bases,
see \cite{nark}, \cite{book}. A number field $K$ of degree $n$ is {\it monogenic} if 
its ring of integers $\Z_K$
is a simple ring extension of $\Z$, that is there exists $\alpha\in\Z_K$ with 
$\Z_K=\Z[\alpha]$. In this case $(1,\alpha,\ldots,\alpha^{n-1})$ is an integral basis of $K$,
called {\it power integral basis}. We also call $\alpha$ the {\it generator}
of this power integral basis.
The algebraic integer $\alpha$ generates a power integral basis
if and only if its {\it index}
\[
I(\alpha)=\sqrt{\left|\frac{D(\alpha)}{D_K}\right|}
\]
is equal to 1, where $D(\alpha)$ is the discriminant of $\alpha$.

The calculation of generators of power integral bases can be reduced to certain
diophantine equations, called {\it index form equations}, cf. \cite{book}.

There exist general algorithms for solving index form equations in cubic, quartic, 
quintic, sextic fields, however the general algorithms for quintic and sextic fields 
are already quite tedious, see \cite{book}.  
Therefore it is worthy to develop efficient methods for
the resolution of special types of higher degree number fields, cf.
\cite{gr}, \cite{book}.

In some cases we considered monogenity in composites of fields, see 
\cite{composite}, \cite{grcompos}. In these cases the index form 
factorizes that makes the resolution of the index form equation easier.

On the other hand, considering totally real relative Thue equations over 
imaginary quadratic fields, it has turned out, that 
the relative Thue equations can be reduced to absolute Thue equations
(over $\Z$), cf. \cite{realthueimquadr}. A similar idea was used in
\cite{totcomplexsextic}, \cite{totcomplexsexticrev}.

In the present paper we study composites $K=LM$ of 
a totally real fields $L$ and imaginary quadratic fields $M$.
We show that in this case the resolution of the index form equation in $K$
can be reduced to solving the index form equation in $L$.
If $L$ is of degree $n$, then the index form equation of $L$ is of degree
$n(n-1)/2$ in $n-1$ variables, 
but the index form equation in $K$ is of degree 
$2n(2n-1)/2$ in $2n-1$ variables. 
Therefore our statement simplifies a lot the calculation
of generators of power integral bases in fields of type $K$ of degree $2n$.

Surprisingly the proofs of our statements are quite simple. However,
they provide a powerful tool.
Our Theorem \ref{main} immediately implies some results of \cite{composite}
(see our remarks after the proof of Theorem \ref{main}).

The strengh and usefulness of our results is aslo demonstrated in the
application given in Section \ref{appl}, where we consider composites of
the so called simplest quartic fields and imaginary quadratic fields.

\section{Composites of totally real fields and imaginary quadratic fields}

Let $f(x)\in\Z[x]$ be a monic irreducible polynomial of degree $n$ having all roots
$\xi=\xi^{(1)},\ldots,\xi^{(n)}$ in $\R$.
The field $L=\Q(\xi)$ is then totally real.
Assume that $L$ has integral basis $(\ell_1=1,\ell_2,\ldots,\ell_n)$ and discriminant
$D_L$.
We shall denote by $\gamma^{(j)}$ the conjugates of any $\gamma\in L \; (j=1,\ldots, n)$.
The index form corresponding to the basis $(\ell_1=1,\ell_2,\ldots,\ell_n)$
is defined by 
\[
I_L(x_2,\ldots,x_n)=\frac{1}{\sqrt{|D_L|}}
\prod_{1\le j_1<j_2\le n}
\left(
x_2\left(\ell_2^{(j_1)}-\ell_2^{(j_2)}\right)+\ldots +x_n\left(\ell_n^{(j_1)}-\ell_n^{(j_2)}\right)
\right).
\]
As it is known, $I_L(x_2,\ldots,x_n)\in\Z[x_2\ldots,x_n]$ and 
for any $x_1,x_2,\ldots,x_n\in\Z$ the algebraic integer element
\[
\alpha=x_1+x_2\ell_2+\ldots x_n\ell_n
\]
generates a power integral basis $(1,\alpha,\ldots,\alpha^{n-1})$ if and only if 
\[
I_L(x_2,\ldots,x_n)=\pm 1.
\]

Let $0<d\in\Z$ be a square-free integer, set 
\[
\omega=
\left\{
\begin{array}{lr}
i\sqrt{d}\;\;& {\rm if}\; -d\equiv 2,3\; (\bmod{4}),\\
\displaystyle{\frac{1+i\sqrt{d}}{2}}\;\; &{\rm if} -d\equiv 1\;  (\bmod{4}).
\end{array}
\right.
\]
The conjugates of $\omega$ are $\omega^{(1)}=\omega,\omega^{(2)}=\overline{\omega}$
(the complex conjugate of $\omega$).
Let $M=\Q(i\sqrt{d})$ with discriminant $D_M$.

Our purpose is to investigate the composite field $K=L\cdot M$ of degree $2\ell$
and discriminant $D_K$.
We assume $(D_L,D_M)=1$. Then
 $(1,\ell_2,\ldots,\ell_n,\omega,\omega\ell_2,\ldots,\omega\ell_n)$
is an integral basis of $K$ and $D_K=D_M^n\cdot D_L^2$ (cf. \cite{nark}).

In the ring of integers $\Z_K$ of $K$ any element $\alpha$ can be 
represented as
\begin{equation}
\alpha=x_1+x_2\ell_2+\ldots x_n\ell_n+y_1\omega+y_2\omega\ell_2+\ldots+y_n\omega\ell_n
=X_1+X_2\ell_2+\ldots+X_n\ell_n,
\label{aa}
\end{equation}
with $x_1,\ldots,x_n,y_1,\ldots,y_n\in\Z$ and $X_j=x_j+\omega y_j$ are elements
in the ring $\Z_M$ of algebraic integers of $M$.

We shall use the following consequence of Theorem 1 of \cite{composite}:

\begin{lemma}
\label{th1}
If $\alpha$ of (\ref{aa}) generates a power integral basis in $K$ then
\begin{equation}
N_{M/Q}(I_L(X_2,\ldots,X_n))=\pm 1
\label{eq1}
\end{equation}
and
\begin{equation}
N_{L/Q}(y_1+y_2\ell_2+\ldots +y_n\ell_n)=\pm 1.
\label{eq2}
\end{equation}
\end{lemma}

The conjugates of $\alpha$ of (\ref{aa}) are obtained obviously as
\[
\alpha^{(j,k)}=x_1+x_2\ell_2^{(j)}+\ldots x_n\ell_n^{(j)}
+y_1\omega^{(k)}+y_2\omega^{(k)}\ell_2^{(j)}+\ldots+y_n\omega^{(k)}\ell_n^{(j)},
\]
for $j=1,\ldots,n,k=1,2$. Note that in this case the index form has three factors, 
two of them are the polynomials on the left hand sides of equations (\ref{eq1}),
(\ref{eq2}). The third factor is 
\begin{equation}
F(x_2,\ldots,x_n,y_1,y_2,\ldots,y_n)=
\prod_{1\le j_1,j_2\le n,j_1\ne j_2} (\alpha^{(j_1,1)}-\alpha^{(j_2,2)}),
\label{third}
\end{equation}
(cf. \cite{composite}, \cite{book}).
This is also a polynomial with coefficients in $\Z$. Equations 
(\ref{eq1}), (\ref{eq2}) together with 
\begin{equation}
F(x_2,\ldots,x_n,y_1,y_2,\ldots,y_n)=\pm 1
\label{eq3}
\end{equation}
are already equivalent to $\alpha$ generating a power integral basis.

In our main result we reduce the relative index form equation (\ref{eq1}) to absolute 
equations, inequalities. This makes the resolution of (\ref{eq1}) much easier.

\begin{theorem}
\label{main}
Assume $\alpha$ of (\ref{aa}) generates a power integral basis in $K$.\\
If $-d\equiv 2,3\;  (\bmod{4})$, then
\begin{equation}
|I_L(x_2,\ldots,x_n)|\le 1
\label{p1}
\end{equation}
and
\begin{equation}
|I_L(y_2,\ldots,y_n)|\le \frac{1}{(\sqrt{d})^{n(n-1)/2}}.
\label{p2}
\end{equation}
If $-d\equiv 1\;  (\bmod{4})$, then
\begin{equation}
|I_L(2x_2+y_2,\ldots,2x_n+y_n)|\le 2^{n(n-1)/2}
\label{p3}
\end{equation}
and
\begin{equation}
|I_L(y_2,\ldots,y_n)|\le \left(\frac{2}{\sqrt{d}}\right)^{n(n-1)/2}.
\label{p4}
\end{equation}
\end{theorem}

\noindent
{\bf Proof of Theorem \ref{main}}. 
According to the arguments in the proof of Theorem 1 of \cite{composite}
we have
\begin{equation}
|N_{M/Q}(I_L(X_2,\ldots,X_n))|=
\prod_{k=1}^2
\left(\frac{1}{|\sqrt{D_L}|}\prod_{1\le j_1<j_2\le n}
\left|\alpha^{(j_1,k)}-\alpha^{(j_2,k)}\right|\right).
\label{ix}
\end{equation}
If $-d\equiv 2,3\;  (\bmod{4})$, then 
\[
{\rm Re}(\alpha^{(j_1,k)}-\alpha^{(j_2,k)})=x_2(\ell_2^{(j_1)}-\ell_2^{(j_2)})+\ldots +
x_n(\ell_n^{(j_1)}-\ell_n^{(j_2)})
\]
and 
\[
{\rm Im}(\alpha^{(j_1,k)}-\alpha^{(j_2,k)})=\sqrt{d}\cdot \left(y_2(\ell_2^{(j_1)}-\ell_2^{(j_2)})+\ldots +
y_n(\ell_n^{(j_1)}-\ell_n^{(j_2)})\right).
\]
We have
\[
|I_L(x_2,\ldots,x_n)|^2
=
|N_{M/Q}(I_L(x_2,\ldots,x_n))|
\]
\[
=\prod_{k=1}^2
\left(\frac{1}{|\sqrt{D_L}|}\prod_{1\le j_1<j_2\le n}
\left|
 x_2(\ell_2^{(j_1)}-\ell_2^{(j_2)})+\ldots + x_n(\ell_n^{(j_1)}-\ell_n^{(j_2)})
\right|\right)
\]
\[
=\prod_{k=1}^2
\left(\frac{1}{|\sqrt{D_L}|}\prod_{1\le j_1<j_2\le n}
\left|{\rm Re}(\alpha^{(j_1,k)}-\alpha^{(j_2,k)})\right|\right)
\]
\[
\le
\prod_{k=1}^2
\left(\frac{1}{|\sqrt{D_L}|}\prod_{1\le j_1<j_2\le n}
\left|\alpha^{(j_1,k)}-\alpha^{(j_2,k)}\right|\right)
\]
\[
= |N_{M/Q}(I_L(X_2,\ldots,X_n))|
\]
whence by (\ref{eq1}) we obtain (\ref{p1}). Similarly,
\[
\left(\sqrt{d}\right)^{n(n-1)}|I_L(y_2,\ldots,y_n)|^2
=
\left(\sqrt{d}\right)^{n(n-1)}|N_{M/Q}(I_L(y_2,\ldots,y_n))|
\]
\[
=\prod_{k=1}^2
\left(\frac{1}{|\sqrt{D_L}|}\prod_{1\le j_1<j_2\le n}
\sqrt{d}\left|
 y_2(\ell_2^{(j_1)}-\ell_2^{(j_2)})+\ldots + y_n(\ell_n^{(j_1)}-\ell_n^{(j_2)})
\right|\right)
\]
\[
=\prod_{k=1}^2
\left(\frac{1}{|\sqrt{D_L}|}\prod_{1\le j_1<j_2\le n}
\left|{\rm Im}(\alpha^{(j_1,k)}-\alpha^{(j_2,k)})\right|\right)
\]
\[
\le
\prod_{k=1}^2
\left(\frac{1}{|\sqrt{D_L}|}\prod_{1\le j_1<j_2\le n}
\left|\alpha^{(j_1,k)}-\alpha^{(j_2,k)}\right|\right)
\]
\[
=  |N_{M/Q}(I_L(X_2,\ldots,X_n))|,
\]
whence by (\ref{eq1}) we obtain (\ref{p2}).

If $-d\equiv 1\;  (\bmod{4})$, then 
\[
{\rm Re}(\alpha^{(j_1,k)}-\alpha^{(j_2,k)})=\frac{2x_2+y_2}{2}\cdot (\ell_2^{(j_1)}-\ell_2^{(j_2)})+\ldots +
\frac{2x_n+y_n}{2}\cdot (\ell_n^{(j_1)}-\ell_n^{(j_2)})
\]
and 
\[
{\rm Im}(\alpha^{(j_1,k)}-\alpha^{(j_2,k)})=\sqrt{d}\cdot \left(\frac{y_2}{2}\cdot (\ell_2^{(j_1)}-\ell_2^{(j_2)})
+\ldots +
\frac{y_n}{2}\cdot (\ell_n^{(j_1)}-\ell_n^{(j_2)})\right).
\]
The above arguments lead us to 
\[
\left|I_L\left(\frac{2x_2+y_2}{2},\ldots,\frac{2x_n+y_n}{2}\right)\right|\le 1
\]
and 
\[
\left|I_L\left(\frac{y_2}{2},\ldots,\frac{y_n}{2}\right)\right|
\le \frac{1}{(\sqrt{d})^{n(n-1)/2}},
\]
whence we obtain (\ref{p3}) and (\ref{p4}) accordingly.
\hfill $\Box$

\vspace{0.5cm}

\noindent
{\bf Remark.}
Our Theorem \ref{main} immediately implies
the result of \cite{composite} on monogenity of 
composites of totally real cyclic fields of prime degree and imaginary quadratic fields, as well as on monogenity of 
composites of Lehmer's quintics and imaginary quadratic fields.
Actually, our result yields, that the above statements of \cite{composite}
are valid for any totally real number field.

\section{Applying Theorem \ref{main}}.

\noindent
Let $-d\equiv 2,3\;  (\bmod{4})$.

If $-d=-1$, then left hand side of both (\ref{p1}) and (\ref{p2}) is $\pm 1$. 
Concerning equation (\ref{p1}) this yields that either $I_(x_2,\ldots,x_n)=0$, that is
$x_2=\ldots=x_n=0$ or $I_(x_2,\ldots,x_n)=\pm 1$, that is $\beta=x_2\ell_2+\ldots+x_n\ell_n$
generates a power integral basis in $L$.
Similary, either $y_2=\ldots=y_n=0$, or $\gamma=y_2\ell_2+\ldots+y_n\ell_n$
generates a power integral basis in $L$. For given $y_2,\ldots,y_n$ we calculate $y_1$ from (\ref{eq2}).
We test all these $x_2,\ldots,x_n,y_1,y_2,\ldots,y_n$ in (\ref{eq3}).
Note that $x_2=\ldots=x_n=0$, and simultaneously $y_2=\ldots=y_n=0$ is not possible, 
since $\omega$ does not generate $K$.

Otherwise, if $-d\neq -1$, (\ref{p2}) gives $y_2=\ldots=y_n=0$ and (\ref{eq2}) gives $y_1=\pm 1$.
The values of $x_2,\ldots,x_n$ are obtained from (\ref{p1}).
If $I_L(x_2,\ldots,x_n)=0$ then $x_2=\ldots x_n=0$ and $\alpha=\pm \omega$ which is again impossible,
since $\omega$ does not generate $K$.
Hence we have to take those $x_2,\ldots,x_n$ for which $I_L(x_2,\ldots,x_n)=\pm 1$, that is 
$\beta=x_2\ell_2+\ldots+x_n\ell_n$ generates a power integral basis in $L$.
We test all these $x_2,\ldots,x_n,y_1,y_2,\ldots,y_n$ in (\ref{eq3}).

\vspace{0.5cm}

\noindent
Let now $-d\equiv 1  (\bmod{4})$. 

If $-d=-3$, then the left hand side of 
(\ref{p3}) is $2^{n(n-1)/2}$, the left hand side of (\ref{p4}) is $(2/\sqrt{3})^{n(n-1)/2}$.
Hence we have to determine the elements $\beta=z_2\ell_2+\ldots+z_n\ell_n$ having index $\leq 2^{n(n-1)/2}$
and to select from those the elements $\gamma=y_2\ell_2+\ldots+y_n\ell_n$ having index 
$\leq (2/\sqrt{3})^{n(n-1)/2}$.  We test if there exist $x_i\in\Z$ with $x_i=(z_i-y_i)/2$.
For given $y_2,\ldots,y_n$ we calculate $y_1$ from (\ref{eq2}).
We test all possible $x_2,\ldots,x_n,y_1,y_2,\ldots,y_n$ in (\ref{eq3}).

Otherwise, if $-d\neq -3$, equation (\ref{p4}) has left hand side $<1$, therefore $y_2=\ldots y_n=0$ and
(\ref{eq2}) gives $y_1=\pm 1$.
(\ref{p3}) implies $|I_L(x_2,\ldots,x_n)|\le 1$. $I_L(x_2,\ldots,x_n)=0$ is not possible again,
$I_L(x_2,\ldots,x_n)=\pm 1$ yields that $\beta=x_2\ell_2+\ldots +x_n\ell_n$ 
generates a power integral basis in $L$.
We test all these $x_2,\ldots,x_n,y_1,y_2,\ldots,y_n$ in (\ref{eq3}).

\vspace{1cm}

This means that for calculating generators of
power integral basis in $K$ we only need the generators of power integral bases in $L$
and, in  case $-d=-3$, elements of small indices of $L$.
This is easy to calculate in lower degree (cubic, quartic) fields, therefore we
obtain an efficient method to calculate generators of power integral bases in sextic, 
octic fields, that are composites of cubic, quartic fields 
and imaginary quadratic fields.

\vspace{0.5cm}

As a consequence of the above arguments we have

\begin{theorem}
\label{cons}
If $-d\ne -1,-3$ then all generators of power integral bases of $K$ are of the form
\[
\alpha=x+\beta\pm \omega
\]
where $\beta$ generates a power integral basis in $L$ and $x\in\Z$ is arbitrary.
\end{theorem}

\vspace{0.5cm}

\noindent
{\bf Remark.} Our Theorem \ref{cons} implies Theorem 2 of \cite{composite} on 
generators of power integral bases in composites of totally real cyclic fields of prime degree and imaginary quadratic fields.

\section{Composites of the simplest quartic fields and imaginary quadratic fields}
\label{appl}

In this section we shall give an application to an infinite parametric family 
of number fields, that shows the strength of our method. 

Let $a>0$ be an integer with $a\ne 3$.  Let $\xi$ be a root of 
\[
f(x)=x^4-ax^3-6x^2+ax+1.
\]
The parametric family of fields $L=\Q(\xi)$ is called {\it simplest quartic fields},
see M. N. Gras \cite{gras}. (She showed that for $a>0,a\ne 3$ the polynomial
$f(x)$ is irreducible.)
In the following we shall also assume that $a^2+16$ is not divisible 
by an odd square. This condition was needed by H. K. Kim and J. H. Lee \cite{kimlee}
to determine an integral basis of $L$. 
(M.N.Gras \cite{gras} showed that $a^2+16$ represents infinitely many
square free integers.)
Using the discriminant 
\[
D(f)=4 (a^2+16)^3
\]
of the polynomial
we obtain the discriminant of $L$. $v_2(x)$ denotes the exponent of 2 in the prime power decomposition of the integer $x$.

\begin{lemma}
\label{kimlemma}
Under the above assumptions an integral basis and the discriminant of $L$ is given by
\[
\begin{array}{rcl}
\left(1,\xi,\xi^2,\frac{1+\xi^3}{2}\right),& D_L=(a^2+16)^3, & {\rm if}\; v_2(t)=0,\\
\left(1,\xi,\frac{1+\xi^2}{2},\frac{\xi+\xi^3}{2}\right),& D_L=\frac{(a^2+16)^3}{4}, &  {\rm if}\; v_2(t)=1,\\
\left(1,\xi,\frac{1+\xi^2}{2},\frac{1+\xi+\xi^2+\xi^3}{4}\right),& D_L=\frac{(a^2+16)^3}{16}, &  {\rm if}\; v_2(t)=2,\\
\left(1,\xi,\frac{1+2\xi-\xi^2}{4},\frac{1+\xi+\xi^2+\xi^3}{4}\right),& D_L=\frac{(a^2+16)^3}{64}, &  {\rm if}\; v_2(t)\geq 3.\\
\end{array}
\]
\end{lemma}

P.Olajos \cite{olaj} determined all generators of power integral bases 
(up to sign  and translation by elements of $\Z$).

\begin{lemma} 
\label{olajlemma}
Under the above assumptions power integral bases exist only for $a=2$ and $a=4$.
All generators of power integral bases are given by
\begin{itemize}
\item
$a = 2$,\;\; $\alpha=x\xi+y\frac{1+\xi^2}{2}+z\frac{\xi+\xi^3}{2}$
where\\
$(x, y, z) = (4, 2,-1), (-13,-9, 4), (-2, 1, 0),
(1, 1, 0),\\ (-8,-3, 2), (-12,-4, 3), (0,-4, 1),
(6, 5,-2), (-1, 1, 0), (0, 1, 0)$
\item
$a = 4$,\;\; $\alpha=x\xi+y\frac{1+\xi^2}{2}+z\frac{1+\xi+\xi^2+\xi^3}{4}$
where\\
$(x, y, z) = (3, 2,-1), (-2,-2, 1),
(4, 8,-3), (-6,-7, 3), \\(0, 3,-1), (1, 3,-1)$.
\end{itemize}
\end{lemma}

Note also that I. Ga\'al and G. Petr\'anyi \cite{gaalpetranyi} calculated the minimal indices
and all elements of minimal index for all parameters $a$.

\vspace{0.5cm}

Our purpose is to study composites of simplest quartic fields with 
imaginary quadratic fields.

Let $d$ be a squarefree positive integer, $M=\Q(i\sqrt{d})$.
Set 
\[
\omega=i\sqrt{d} \;\; {\rm if} \;\; -d\equiv 2,3\; (\bmod{4})\;\;
{\rm and}\;\;
\omega=\frac{1+i\sqrt{d}}{4} \;\; {\rm if} \;\; -d\equiv 1\; (\bmod{4}).
\]
We assume that $(D_M,D_L)=1$. If $v_2(a)\geq 1$ then $D_L$ is even, hence
in that case we must have $-d\equiv 1\; (\bmod{4})$.

We consider monogenity of the composite field $K=LM$ of degree 8.
This is an infinite parametric family of octic fields, depending 
on the parameters $d,a$.

We show:

\begin{theorem}
\label{cq}
Under the above conditions for $d\neq 3$ the field $K$ is not monogenic.
\end{theorem}

\noindent
{\bf Proof of Theorem \ref{cq}}

For a given parameter $a$ denote by $(1,\ell_2,\ell_3,\ell_4)$ 
an integral basis of $L$. Denote by $I_L(x_2,x_3,x_4)$ the index form
corresponding to this integral basis.
The condition $(D_L,D_M)=1$ implies that an integral basis of $K$ is
given by
$(1,\ell_2,\ell_3,\ell_4,\omega,\omega\ell_2,\omega\ell_3,\omega\ell_4)$.
Represent any $\alpha\in\Z_K$
in the form
\[
\alpha=x_1+x_2\ell_2+x_3\ell_3+x_4\ell_4
+y_1\omega +y_2\omega\ell_2+y_3\omega\ell_3+y_4\omega\ell_4
\]
with integer coefficients $x_i,y_i\; (1\leq i\leq 4)$.

{\bf I.} Assume $d\neq 1$. By Theorem \ref{main}
we have (in both cases $-d\equiv 2,3\; (\bmod{4})$
and $-d\equiv 1\; (\bmod{4})$) 
$I_L(y_2,y_3,y_4)=0$ implying $y_2=y_3=y_4=0$. Theorem \ref{th1},
equation (\ref{eq2}) implies $y_1=\pm 1$. Further, by Theorem \ref{main}
$|I_L(x_2,x_3,x_4)|\leq 1$. 
$I_L(x_2,x_3,x_4)=0$ would imply $x_2=x_3=x_4=0$ which is impossible,
while $\omega$ does not generate $K$. Therefore 
$I_L(x_2,x_3,x_4)=\pm 1$ which means that $\beta=x_2\ell_2+x_3\ell_3+x_4\ell_4$
generates a power integral basis in $L$. 
By Lemma \ref{olajlemma} this is only possible for $a=2,4$ in which cases
all possible $x_2,x_3,x_4$ are known.

Consider now the third factor (\ref{eq3}) of the index form equation of $K$.
In both cases $-d\equiv 2,3\; (\bmod{4})$ and $-d\equiv 1\; (\bmod{4})$,
both for $a=2$ and $a=4$ we substitute the possible $x_2,x_3,x_4$ and
$y_2=y_3=y_4=0$, $y_1=\pm 1$ into $F(x_2,x_3,x_4,y_1,y_2,y_3,y_4)$.
In all cases we obtain a polynomial of degree 6 in $d$ with positive
coefficients, implying that $F(x_2,x_3,x_4,y_1,y_2,y_3,y_4)=\pm 1$
is not possible.

{\bf II.} Let now $d=1$. 
By Theorem \ref{main} we have
\[
|I_L(x_2,x_3,x_4)|\leq 1,\;\;\; |I_L(y_2,y_3,y_4)|\leq 1.
\]

\noindent
a.) $I_L(x_2,x_3,x_4)=I_L(y_2,y_3,y_4)=0$ implies $x_2=x_3=x_4=0$ and 
$y_2=y_3=y_4=0$, this is not possible, while $\omega$ does not generate $K$.
Hence we must have either $I_L(x_2,x_3,x_4)=\pm 1$ or $I_L(y_2,y_3,y_4)=\pm 1$
(or both), hence there must be a power integral basis in $L$, whence
only $a=2$ and $a=4$ is possible. In the following we restrict ourselves
to $a=2$ and $a=4$.

\noindent
b.) Assume $I_L(x_2,x_3,x_4)=\pm 1$ and $I_L(y_2,y_3,y_4)=0$.
The possible $x_2,x_3,x_4$ are listed in Lemma \ref{olajlemma},
$y_2=y_3=y_4=0$, $y_1=\pm 1$. We substitute the possible values of 
the variables into $F(x_2,x_3,x_4,y_1,y_2,y_3,y_4)$ and 
find that in all cases it takes huge values, not $\pm 1$.

\noindent
c.) Assume $I_L(x_2,x_3,x_4)=0$ and $I_L(y_2,y_3,y_4)=\pm 1$.
Then $x_2=x_3=x_4=0$ and the possible $y_2,y_3,y_4$ are listed in 
Lemma \ref{olajlemma}. 
We substitute the possible values $x_2=x_3=x_4=0$ and $y_2,y_3,y_4$
into $F(x_2,x_3,x_4,y_1,y_2,y_3,y_4)$
and obtain a polynomial of degree 12 in $y_1$. 
Solving $F(x_2,x_3,x_4,y_1,y_2,y_3,y_4)=\pm 1$ in $y_1$ 
we never get integer solutions for $y_1$.

\noindent
d.) If $I_L(x_2,x_3,x_4)=\pm 1$ and $I_L(y_2,y_3,y_4)=\pm 1$,
then $x_2,x_3,x_4$ and $y_2,y_3,y_4$ may run independently through
the possible triplets listed in Lemma \ref{olajlemma}.
In all cases we substitute $x_2,x_3,x_4$ and $y_2,y_3,y_4$
into $F(x_2,x_3,x_4,y_1,y_2,y_3,y_4)$
and obtain a polynomial of degree 12 in $y_1$. 
Solving $F(x_2,x_3,x_4,y_1,y_2,y_3,y_4)=\pm 1$ in $y_1$ 
we never get integer solutions for $y_1$.
\hfill $\Box$

\vspace{0.5cm}

\noindent
{\bf Remark.} The case $a=3$ is not covered by Theorem \ref{cq}. 
In this case we have
\[
|I_L(2x_2+y_2,2x_3+y_3,2x_4+y_4)|\leq 64
\]
and 
\[
|I_L(y_2,y_3,y_4)|\leq 2.
\]
According to \cite{gaalpetranyi} we have minimal index 2 in $L$ 
for infinitely many parameters $a$. Especially the elements
of index $\leq 64$ seems very difficult to determine. 
This could be the subject of a further research.

\vspace{0.5cm}

\noindent
{\bf Remark.} Note that in \cite{grcompos} we obtained conditions on the monogenity of 
composites of fields, among others of simplest quartic fields and quadratic fields. 
We did not assume that the discriminants are relative prime and involved also real
quadratic fields. However we only obtained certain divisibility conditions on
the parameters as a consequence of monogenity.

\section{One more example}

We also provide a positive example to show that such composite fields
may happen to be monogenic. Consider the totally real quartic field $L$ generated by a root
$\xi$ of the polynomial $f(x)=x^4-4x^2-x+1$. In this field $(1,\xi,\xi^2,\xi^2)$ is an integral basis,
$D_L=1957$. Let $M=\Q(i)$ with $D_M=-4$, coprime to $D_L$.
The composite field $K=LM$ can be generated e.g. by $\alpha=i\xi$ having minimal polynomial 
$g(x)=x^8+8x^6+18x^4+9x^2+1$. In $K$ the element $\alpha$ generates a power integral basis.

\section{Computational aspects}

All calculations connected with the above examples
were performed in Maple \cite{maple}.
Our procedures were executed 
on an average laptop running under Windows. The CPU time took 
all together some seconds.

\end{document}